\documentclass[12pt,reqno]{amsart}

\usepackage{lipsum}
\usepackage{graphicx}
\usepackage{array}
\usepackage{amssymb}
\usepackage{hyperref}
\usepackage{amsthm}
\usepackage{amsfonts}
\usepackage{amsmath}
\usepackage{bm}
\usepackage{mathrsfs}
\usepackage[all]{xy}
\usepackage{color}
\usepackage{subfigure}
\usepackage{tikz, tikz-cd}
\usepackage{enumerate}
\usepackage{anysize}
\usepackage{amscd}
\usepackage[letterpaper]{geometry}
\usepackage{geometry}
\usepackage{multirow}
\usepackage{array}
\usepackage{booktabs}
\usepackage{musicography}
\usepackage[utf8]{inputenc}
\usepackage[T1]{fontenc}

\newtheorem{theorem}{Theorem}
\newtheorem{remark}{Remark}

\newtheorem{lemma}[theorem]{Lemma}
\newtheorem{proposition}[theorem]{Proposition}
\theoremstyle{definition}

\allowdisplaybreaks[3]

\begin{document}

 \title{Spin(7)-Manifolds and Multisymplectic Geometry}

 \author{Aaron Kennon} \thanks{The research of the author is partially suported by Simons Foundation Award \#488629 (Morrison)}
  \address{Department of Physics UCSB}
 \email{akennon@physics.ucsb.edu}

\begin{abstract}
We utilize Spin(7) identities to prove that the Cayley four-form associated to a torsion-free Spin(7)-Structure is non-degenerate in the sense of multisymplectic geometry. Therefore, Spin(7) geometry may be treated as a special case of multisymplectic geometry. We then capitalize on this relationship to make statements about Hamiltonian multivector fields and differential forms associated to torsion-free Spin(7)-Structures. 
\end{abstract}

\date{}
\maketitle
\setcounter{tocdepth}{1}

\section{Introduction}

The Berger classification states that the possible holonomy groups for simply-connected, irreducible, non-symmetric Riemannian manifolds are quite constrained \cite{Berger}. According to this classification, Spin(7) is one such possible holonomy group for eight-dimensional Riemannian manifolds.  Such manifolds are also Ricci-Flat and admit precisely one parallel spinor, making them interesting backgrounds for M-Theory and F-Theory \cite{GukovSparks, BGPP}. \\

Despite the significant role of these manifolds in both differential geometry and physics we know relatively little about them. Much of what we do know comes from generic results in holonomy theory, spin geometry, and from what we can infer from specific constructions. Bryant proved local existence of Spin(7)-holonomy metrics \cite{Bryant} and then with Salamon constructed the first complete non-compact example \cite{BryantSalamon}. They constructed a one-parameter family of complete Spin(7)-metrics on the total space of the spinor bundle over S$^{4}$ and these examples have asymptotically conical geometry.  The first compact examples were constructed by Joyce as resolutions of toroidal orbifolds \cite{Joyce1}. Since then many more examples have been constructed, including non-compact examples with asymptotically locally conical geometry \cite{CGLP, CGLP2, GukovSparks, Foscolo} and more compact examples by Joyce utilizing certain Calabi-Yau four-orbifolds with antiholomorphic in\textnormal{vol}utions \cite{Joyce2}. \\

As exceptional holonomy manifolds are associated to closed differential forms satisfying non-denegeneracy conditions there is at least a tentative analogy between exceptional holonomy and symplectic geometry. In comparision a lot more is known about symplectic geometry than exceptional holonomy manifolds, so one way of potentially understanding these exceptional holonomy manifolds better is to see how far we can push this analogy. As the differential forms associated to exceptional holonomy manifolds naturally define multisymplectic structures, studying exceptional holonomy manifolds in this context may clarify the relationships between these spaces and symplectic manifolds. This programme to approach exceptional geometry from the perspective of multisymplectic geometry has been explored in the G$_{2}$-holonomy context (and that of closed G$_{2}$-Structures and co-closed G$_{2}$-Structures) by Cho, Salur, and Todd \cite{CST}. In this work we investigate if similar results hold in the Spin(7) context. \\

Recall that in symplectic geometry we have a closed two-form $\omega$ which satisfies a nondegeneracy condition. 

\begin{equation*}
X \lrcorner \omega = 0 \hspace{5mm}  \textnormal{iff} \hspace{5mm}  X=0
\end{equation*}

On a manifold which admits a symplectic form $\omega$, a vector field $X_{H}$ and a function $H$ are said to be Hamiltonian if the vector field contracted into the symplectic form gives the differential of the function. 

\begin{equation*}
X_{H} \lrcorner \omega = dH
\end{equation*}

Like symplectic structures, multisymplectic structures are characterized by a closed differential form satisfying a nondegeneracy condition, however, we allow the form associated to a multisymplectic structure to have arbitrary degree. Given such a multisymplectic structure we also have analogous notions of Hamiltonian multivector fields and differential forms. Just as symplectic geometry provides the natural model for phase space of Hamiltonian systems, multisymplectic geometry comes up naturally in physics in the covariant Hamiltonian formulation of classical mechanics \cite{CIdL1, CIdL2, CCI, FPR, PR}. Multisymplectic geometry as a discipline in and of itself is therefore highly relevant for studying quantization of classical field theories. Key results from symplectic geometry such as the Darboux theorem and momentum mapping construction in certain instances generalize to multisymplectic geometry, although the results are usually not quite as strong as in the symplectic case as multisymplectic geometry is a much more general framework \cite{RW}.  \\ 

Given an eight-manifold with a torsion-free Spin(7)-Structure, which is a more general structure than a Spin(7) holonomy metric, we have a closed four-form called the Cayley form which turns out to satisfy the multisymplectic nondegeneracy condition. This motivates us to study the Hamiltonian multivector fields and differential forms associated to the Cayley four-form, which we call Cayley multivector fields and differential forms to distinguish them from the more general context of multisymplectic geometry. \\

We construct maps $\Psi_{k}$ for k equals 1, 2, and 3 by contracting k-multivector fields with the Cayley four-form. We prove that $\Psi_{1}$ is injective, $\Psi_{2}$ is an isomorphism, and that $\Psi_{3}$ is surjective. Since the latter two maps are (at least) surjective, it follows that Cayley two-multivector fields and three-multivector fields exist on a manifold with a torsion-free Spin(7)-Structure. However, we prove a general nonexistence result for Cayley vector fields on closed Spin(7)-Manifolds. We then characterize the decomposition of Cayley two-multivector fields and three-multivector fields into irreducible Spin(7) representations and prove some results constraining these multivector fields and differential forms.  \\

We discuss background needed for Spin(7) geometry and multisymplectic geometry in \S2 and our results in \S3. 

\section{Background}

\subsection{Spin(7) Geometry}

We begin with an introduction to Spin(7) geometry. The material presented here and more general background may be found in Joyce's book \cite{JoyceBook} and Karigiannis' thesis \cite{Karigiannis}. \\

The group Spin(7) and more specifically its representations plays a central role in the geometry of these structures. Spin(7) may be defined as the double cover of SO(7), but for our purposes it is most relevant that it is the automorphism group of the four-form $\Psi_{0}$ on $\mathbb{R}^{8}$. 

\begin{equation*}
\begin{aligned}
\Psi_{0} &= dx^{0123} - dx^{0167} - dx^{0527} - dx^{0563} - dx^{0415} - dx^{0426} - dx^{0437} \\
      &+ dx^{4567} - dx^{4523} - dx^{4163} - dx^{4127} - dx^{2637} - dx^{1537} - dx^{1526}
\end{aligned}
\end{equation*}

Here, the notation $dx^{ijkl}$ is shorthand for $dx^{i} \wedge dx^{j} \wedge dx^{k} \wedge dx^{l}$. This four-form is referred to as the canonical Spin(7)-Structure on $\mathbb{R}^{8}$. It important to note that this form is self-dual and it defines a Riemannian metric. This metric is flat and as a consequence this four-form is also sometimes called the flat Spin(7) form. \\

The group Spin(7) is relevant to us because it is one of the possible holonomy groups for an eight-dimensional irreducible Riemannian manifold with a parallel spinor. The generic orientable Riemannian eight-manifold has holonomy SO(8), but we can ensure a holonomy reduction by showing the manifold admits a torsion-free Spin(7)-Structure. Existence of a Spin(7)-Structure more generally is just a topological condition, it is equivalent to the manifold being orientable, spinnable, and satisfying an additional characteristic class condition.

\begin{equation*}
p_{1}(M) -4p_{2}(M) \pm 8 \chi (M) = 0
\end{equation*}

We can also characterize a Spin(7)-Structure in terms of a specific four-form that is pointwise identified with the canonical four-form in $\mathbb{R}^{8}$. Such a four-form is called \textit{admissible}. Since the canonical four-form is self-dual, this implies the same for any admissible four-form. Note that unlike the analogous notion of positivity for a G$_{2}$-Structure, admissability is not an open condition. We call an admissible four-form a \textit{Cayley four-form} to emphasize the relationship of these structures to the octonions. \\

The existence of an admissible four-form by itself is not very special, but it turns out that if the admissible four-form is closed it corresponds to a torsion-free Spin(7)-Structure. An eight-dimensional manifold with a torsion-free Spin(7)-Structure will be called a \textit{Spin(7)-Manifold}. From the definition it is possible that a Spin(7)-Manifold may have holonomy properly contained in Spin(7). It is often important in applications to ensure that the holonomy is precisely Spin(7) and we can do so by checking some additional criteria. In the simply-connected compact case we can show that if the holonomy is known to be contained in Spin(7), the holonomy is precisely Spin(7) if and only if the $\hat{A}$ genus is one. We can then use results from spin geometry to rephrase this constraint on the $\hat{A}$ genus in terms of a relationship between the Betti numbers of the manifold.

\begin{equation*}
25=-b_{2}+b_{3}+b^{+}_{4}-2b^{-}_{4}
\end{equation*}

In the compact simply-connected case, the Betti numbers must be related in this manner for a manifold to admit a Spin(7)-holonomy metric. It is also known that compact manifolds with Spin(7)-holonomy metrics are simply connected. For a more general torsion-free Spin(7)-Structure we also have a constraint relating the $\hat{A}$ genus to the Betti numbers and in this case the first Betti number doesn't vanish. \\

Given a Spin(7)-Structure, not necessarily torsion-free, we can decompose the spaces of k-forms orthogonally into irreducible representations of Spin(7). If $\Lambda^{k}(T^{*}M)$ denotes the space of k-forms we use the notation $\Lambda^{k}_{\ell}(T^{*}M)$ to denote the representation space of dimension $\ell$ that is a summand of $\Lambda^{k}(T^{*}M)$. The spaces of one-forms and seven-forms each consist of one irreducible Spin(7) representation whereas for intermediary degree forms we get a direct sum of irreducible representations. 

\begin{align*}
 \hspace{-100mm} & \Lambda^{1}(T^{*}M) = \Lambda^{1}_{8}(T^{*}M) \\
& \Lambda^{2}(T^{*}M) = \Lambda^{2}_{7}(T^{*}M) \oplus \Lambda^{2}_{21}(T^{*}M) \\
& \Lambda^{3}(T^{*}M) = \Lambda^{3}_{8}(T^{*}M) \oplus \Lambda^{3}_{48}(T^{*}M) \\
& \Lambda^{4}(T^{*}M) = \Lambda^{4}_{1}(T^{*}M) \oplus \Lambda^{4}_{7}(T^{*}M) \oplus \Lambda^{4}_{27}(T^{*}M) \oplus \Lambda^{4}_{35}(T^{*}M)  \hspace{23mm}
\end{align*}

The representations $\Lambda^{k}(T^{*}M)$ for k equals 5, 6, 7 and 8 decompose in precisely the same manner as the corresponding $\Lambda^{8-k}(T^{*}M)$ and the representations are related by the Hodge star. We can describe these representation spaces explicitly in terms of the admissible four-form $\Psi$ associated to the Spin(7)-Structure \cite{Karigiannis}.

\begin{align*}
& \Lambda^{2}_{7}(T^{*}M) = \{ \beta \in \Lambda^{2}(T^{*}M) : *(\Psi \wedge \beta) = -3\beta \} \\
& \Lambda^{2}_{21}(T^{*}M) =  \{ \beta \in \Lambda^{2}(T^{*}M) : *(\Psi \wedge \beta) = \beta \} \\ 
& \Lambda^{3}_{8}(T^{*}M)= \{ w \lrcorner \Psi \hspace{2mm} \textnormal{for w} \in \Gamma(TM) \} \\
& \Lambda^{3}_{48}(T^{*}M) =  \{ \eta \in \Lambda^{3}(T^{*}M): \eta \wedge \Psi=0 \} \\ 
& \Lambda^{4}_{1}(T^{*}M) = \{ f\Psi \hspace{2mm}\textnormal{for} \hspace{2mm} f \in C^{\infty}(M) \} \\
& \Lambda^{4}_{7}(T^{*}M) =  \{ v^{\#}\wedge (w \lrcorner \Psi)- w^{\#} \wedge (v \lrcorner \Psi) \hspace{2mm} \textnormal{for v,w} \in \Gamma(TM) \} \\
& \Lambda^{4}_{27}(T^{*}M) =  \{ \sigma \in \Lambda^{4}(T^{*}M) : *\sigma=\sigma, \sigma \wedge \Psi=0, \sigma \wedge \tau=0 \hspace{5mm} \forall \tau \in \Lambda^{4}_{7}(T^{*}M)\} \\
& \Lambda^{4}_{35}(T^{*}M) =  \{ \sigma \in \Lambda^{4}(T^{*}M) : *\sigma=-\sigma \}
\end{align*}

In certain cases we can explicitly write down isomorphisms of appropriate representations by wedging the differential forms with $\Psi$. For instance, the representation space $ \Lambda^{4}_{1}(T^{*}M)$ is isomorphic to the space of functions on the underlying manifold and $\Lambda^{5}_{8}(T^{*}M)$ is isomorphic to the space of one-forms. \\

\subsection{Multisymplectic Geometry}

As we have seen, we may associate to a torsion-free Spin(7)-Structure a closed, admissible four-form. The fact that this four-form has these properties allows us to make contact with multisymplectic geometry. Useful references for multisymplectic geometry material here and more generally include the recent review by Ryvkin and Wurzbacher \cite{RW} and the treatment from Cantrijin et al \cite{CIdL1, CIdL2}. \\

We call an n-dimensional manifold equipped with a k+1 degree closed differential form $\xi$ a \textit{k-multisymplectic n-manifold} if the interior product with $\xi$ as a map is injective. 

\begin{equation*}
X \lrcorner \xi = 0 \hspace{5mm}  \textnormal{iff} \hspace{5mm}  X=0
\end{equation*}

This condition is called a nondegeneracy condition and it is precisely the same condition as that on the two-form in symplectic geometry. In this language a symplectic structure is a one-multisymplectic structure and the existence of such a structure implies that the underlying manifold is even dimensional. It turns out that admissibility of a Spin(7)-Structure implies this notion of nondegeneracy used in multisymplectic geometry. This is demonstrated using the Spin(7) decomposition of differential forms in Lemma 1 of \S3. \\

Provided that the Spin(7)-Structure is torsion-free, the Cayley four-form is also closed so in this case a torsion-free Spin(7)-Structure defines a three-multisymplectic structure on an eight-dimensional manifold. In the G$_{2}$ case certain lower torsion classes also define multisymplectic structures, not just the torsion-free class. In particular, a closed G$_{2}$-Structure defines a two-multisymplectic structure and a co-closed G$_{2}$-Structure defines a three-multisymplectic structure. Analogously to how admissibility in the Spin(7) case generally implies multisymplectic nondegeneracy, G$_{2}$ positivity also implies nondegeneracy. \\

Multisymplectic Geometry, as its name suggests, features certain tensors called multivector fields. The construction of the vector space of k-multivector fields is analogous to that of the space of differential k-forms. Recall that the space of differential k-forms is by definition the space of sections of the k$^{\textnormal{th}}$ exterior power of the cotangent bundle. We construct the space of k-multivector fields using a similar procedure, but instead of starting with the cotangent bundle we take the k$^{\textnormal{th}}$ exterior power of the tangent bundle. We can put the spaces of k-multivector fields or differential k-forms of various degrees together into associative graded algebras where the product operation is the wedge product. 

\begin{align*}
\Lambda(TM)&=\bigoplus_{k=0}^{n} \Lambda^{k}(TM) \\
\Lambda(T^{*}M)&=\bigoplus_{k=0}^{n} \Lambda^{k}(T^{*}M)
\end{align*}

We say a k-multivector field $Q$ is \textit{decomposable} if it may be written as the wedge product of k individual vector fields.

\begin{equation*}
Q= u_{1} \wedge ...\wedge u_{k}
\end{equation*}

Decomposable multivector fields play an important role in the theory as any multivector field may be written as a sum of decomposable terms. This fact is used to define the interior product of an $\ell$-multivector field with a k-form. For $Q = u_{1} \wedge ...\wedge u_{\ell}$ a decomposable l-form we may define the interior product of $Q$ with a k-form $\beta$ in terms of the usual interior product of vector fields \cite{CIdL1}. 

\begin{equation*}
Q \lrcorner \beta=u_{\ell}\lrcorner...u_{1} \lrcorner \beta
\end{equation*}

The resulting differential form is a (k-$\ell$)-form. Since all multivector fields are built out of decomposable parts, we may define this contraction operation for an arbitrary multivector field by extending by linearity. \\

Given a Riemannian n-manifold there are a variety of identities relating the interior product of a vector field $X$ into a differential k-form $\beta$ to the Hodge star and metric dual operations \cite{Karigiannis}. 

\begin{align*}
& X \lrcorner \beta=(-1)^{(n-k)(k-1)} \ast (X^{\musFlat{}} \wedge \ast\beta) \\
& X \lrcorner \ast \beta=(-1)^{k} \ast (X^{\musFlat{}} \wedge \beta) \\
& \ast(X \lrcorner \beta)=(-1)^{(k+1)} X^{\musFlat{}} \wedge \ast\beta \\
& \ast(X \lrcorner \ast\beta)=(-1)^{n-k-1+k(n-k)} X^{\musFlat{}} \wedge \beta 
\end{align*}

The Hodge star and musical isomorphisms extend immediately to the multivector field context which allows us to prove similiar identities on an n-dimensional manifold for an $\ell$-multivector field Q and k-form $\beta$ \cite{CST}.

\begin{align*}
& Q \lrcorner \beta=(-1)^{(k-l)(n-k)} \ast (Q^{\musFlat{}} \wedge \ast\beta) \\
& Q \lrcorner \ast \beta=(-1)^{kl} \ast (Q^{\musFlat{}} \wedge \beta) \\
& \ast(Q \lrcorner \beta)=(-1)^{l(k-l)} Q^{\musFlat{}} \wedge \ast\beta \\
& \ast(Q \lrcorner \ast\beta)=(-1)^{l(n-k-l)+k(n-k)} Q^{\musFlat{}} \wedge \beta 
\end{align*}

The central objects in multisymplectic geometry are Hamiltonian multivector fields and differential forms. For k-multisymplectic structure $\xi$, a (k-$\ell$)-multivector field $Q$ is called \textit{locally Hamiltonian} if the form $Q \lrcorner \xi$ is closed. $Q$ is called a \textit{Hamiltonian (k-$\ell$)-multivector field} if for some $\alpha$ $\in$ $\Lambda^{\ell}(T^{*}M)$ the contraction of $Q$ with $\xi$ gives an exact form. 

\begin{equation*}
Q \lrcorner \xi = d\alpha
\end{equation*}

Likewise if for $\alpha$ $\in$ $\Lambda^{\ell}(T^{*}M)$ there exists a $Q$ $\in$ $\Lambda^{k-\ell}(TM)$ that satisfies the above equation, then $\alpha$ is called a \textit{Hamiltonian differential $\ell$-form}. \\

Note that a Hamiltonian multivector field doesn't determine the corresponding Hamiltonian differential form uniquely, it only does so up to the addition of a closed form.  Similiarly, a Hamiltonian differential form only determines the corresponding multivector field up to the addition of a multivector field whose contraction with the k-multisymplectic form is zero. We mention for the sake of completeness that one can remove this nonuniqueness by quotienting out from the space of Hamiltonian differential k-forms the space of closed k-forms and from the space of Hamiltonian k-multivector fields the space of multivector fields whose contraction with the multisymplectic form vanishes. The resulting quotient spaces will then be isomorphic.  \\

We may also extend other operations associated with vectors to the multivector context. The generalization of the Lie bracket of vector fields to multisymplectic geometry is called the \textit{Schouten-Nijenhuis bracket}. We first define it when the multivector field associated with the Lie derivative is decomposable. In that case we define the Lie bracket of a multivector field $Q$ with respect to a decomposable multivector field in terms of the Lie derivatives of $Q$ with respect to the individual vector fields. 

\begin{equation*}
[X_1 \wedge \cdots \wedge X_{l},Q]=\sum_{i=0}^{l}(-1)^{i}X_1\wedge \cdots \wedge \hat{X}_i\wedge \cdots \wedge X_{l}\wedge L_{X_{i}} Q
\end{equation*}

We can then extend this bracket by linearity to arrive at the Schouten-Nijenhus bracket for arbitrary multivector fields. The Schouten-Nijenhuis bracket has several properties reminiscent of the Lie bracket \cite{CIdL1}. 

\begin{equation*}
[Q_1,Q_2]=(-1)^{q_1q_2}[Q_2,Q_1]
\end{equation*}
\begin{equation*}
[Q_1,Q_2\wedge Q_3]=[Q_1,Q_2]\wedge Q_3 + (-1)^{q_1q_2+q_2}Q_2\wedge[Q_1,Q_3]
\end{equation*}
\begin{equation*}
\begin{split}
(-1)^{q_1(q_3-1)}&[Q_1,[Q_2,Q_3]]+(-1)^{q_2(q_1-1)}[Q_2,[Q_3,Q_1]] \\
&+(-1)^{q_3(q_2-1)}[Q_3,[Q_1,Q_2]]=0 \\
\end{split}
\end{equation*}

We can also extend the concept of the Lie derivative more generally to the multivector context. For instance we define the Lie derivative of a differential form $\beta$ with respect to an $\ell$-multivector field $Q$ using an extension of Cartan's formula for the standard Lie derivative acting on differential forms. 

\begin{equation*}
 L_Q\beta=Q \lrcorner d\beta - (-1)^q d(Q \lrcorner \beta)
\end{equation*}

This Lie derivative satisfies many useful properties relating it to the exterior derivative and Schouten-Nijenhuis bracket \cite{CIdL1}. 

\begin{equation*}
 d L_Q\beta=(-1)^{q+1} L_Q d\beta
\end{equation*}
\begin{equation*}
 L_{Q_1\wedge Q_2}\beta=Q_2 \lrcorner L_{Q_1}\beta+(-1)^{q_1}L_{Q_2}Q_1 \lrcorner\beta
\end{equation*}
\begin{equation*}
 [Q_1,Q_2] \lrcorner\beta=(-1)^{q_1q_2+q_2} L_{Q_1}Q_2\lrcorner\beta-Q_2\lrcorner L_{Q_1}\beta
\end{equation*}

We will see that the Lie derivative with respect to a Hamiltonian multivector field will always vanish when acting on the multisymplectic form. This fact turns out to be particularly relevant in the case of a Hamiltonian vector field corresponding to a torsion-free Spin(7)-Structure. 

\section{Results}

We now specialize to the case that the multisymplectic form is the Cayley four-form associated to a torsion-free Spin(7)-Structure. In this case since the relevant differential form is a four-form we have three associated linear maps given by contracting the four-form with vector fields, two-multivector fields, and three-multivector fields.

\begin{align*}
& \Psi_{1}: \Gamma(TM) \rightarrow \Lambda^{3}(T^{*}M) \\
& \Psi_{2}: \Lambda^{2}(TM) \rightarrow \Lambda^{2}(T^{*}M) \\
& \Psi_{3}: \Lambda^{3}(TM) \rightarrow \Lambda^{1}(T^{*}M) 
\end{align*}

We first prove that the map $\Psi_{1}$ is injective, which implies that a torsion-free Cayley four-form is a three-multisymplectic structure.  

\begin{lemma}
The map $\Psi_{1}: \Gamma(TM) \rightarrow \Lambda^{3}(T^{*}M)$ is injective.
\end{lemma}

\textit{Proof.} We start with the standard formula relating the interior product to the Hodge star and musical isomorphism. 

\begin{equation*}
X \lrcorner \Psi = * (\Psi \wedge X^{\musFlat{}})
\end{equation*}

Therefore, if $X \lrcorner \Psi$ vanishes we may infer that the one-form wedged with the Cayley form vanishes since the Hodge star is a linear isomorphism. \\

\begin{equation*}
\Psi \wedge X^{\musFlat{}} = 0
\end{equation*}

However, the map from $\Lambda^{1}_{8}(T^{*}M)$ to $\Lambda^{5}_{8}(T^{*}M)$ given by wedging the one-form with $\Psi$ is an isomorphism of representations, so the vanishing of $\Psi \wedge X^{\musFlat{}}$ implies that $X^{\musFlat{}} = X = 0$. Therefore, if a vector field contracted with the Cayley form vanishes then that vector field must be identically zero. This condition is equivalent to the map $\Psi_{1}$ being injective. \\

\begin{lemma}
The map $\Psi_{3}: \Lambda^{3}(TM) \rightarrow \Lambda^{1}(T^{*}M) $ is surjective.
\end{lemma}

\textit{Proof.} We can express any given one-form $\alpha$ in terms of itself and the Cayley four-form in a con\textnormal{vol}uted manner \cite{Karigiannis}.

\begin{equation*}
\ast (\Psi \wedge \ast(\Psi \wedge \alpha))=-7\alpha
\end{equation*}

We may then write $\ast (\Psi \wedge \alpha)$ in terms of the standard interior product

\begin{equation*}
\ast (\Psi \wedge \alpha) = \alpha^{\musSharp{}} \lrcorner \Psi
\end{equation*}

$\alpha^{\musSharp{}} \lrcorner \Psi$ is a three-form and its metric dual is a three-multivector field we denote as Q. We can then rewrite the $\alpha$ identity in terms of Q. 

\begin{equation*}
\ast (\Psi \wedge Q^{\musFlat{}}) =-7\alpha
\end{equation*}

This expression is in a form that we can rewrite in terms of the multisympectic interior product.

\begin{equation*}
\big(-\frac{1}{7}Q\big) \lrcorner \Psi = \alpha
\end{equation*}

Therefore, the map $\Psi_{3}$ is surjective as any one-form $\alpha$ may be written in this manner. \\

Now we may prove that the map $\Psi_{2}$ has the best qualities of $\Psi_{1}$ and $\Psi_{3}$. \\

\begin{lemma}
The map $\Psi_{2}: \Lambda^{2}(TM) \rightarrow \Lambda^{2}(T^{*}M)$ is an isomorphism.
\end{lemma}

\textit{Proof.} We will first show the map is injective and to do so we will look at the kernel which consists of the two-multivector fields $Q$ whose contraction with $\Psi$ vanishes. 

\begin{equation*}
Q \lrcorner \Psi = 0
\end{equation*}

As we have seen, we can rewrite the above interior product equation in terms of an equation in\textnormal{vol}ving the Hodge star and metric dual to $Q$.

\begin{equation*}
\ast (\Psi \wedge Q^{\musFlat{}}) = Q \lrcorner \Psi = 0
\end{equation*}

We can then decompose the two-form $Q^{\musFlat{}}$ into irreducible Spin(7) representations.

\begin{equation*}
Q^{\musFlat{}} = Q^{\musFlat{}}_{7}+Q^{\musFlat{}}_{21}
\end{equation*}

$Q^{\musFlat{}}_{7}$ and $Q^{\musFlat{}}_{21}$ satisfy the equations which define their representations. 

\begin{align*}
& \ast (\Psi \wedge Q^{\musFlat{}}_{7})=-3Q^{\musFlat{}}_{7} \\
& \ast (\Psi \wedge Q^{\musFlat{}}_{21})=Q^{\musFlat{}}_{21}
\end{align*}

By assumption that $Q$ is in the kernel of $\Psi_{2}$ we see that a linear combination of $Q^{\musFlat{}}_{7}$ and $Q^{\musFlat{}}_{21}$ has to vanish. 

\begin{equation*}
0 = Q \lrcorner \Psi = \ast (\Psi \wedge Q^{\musFlat{}}) = \ast (\Psi \wedge Q^{\musFlat{}}_{7}) + \ast (\Psi \wedge Q^{\musFlat{}}_{21})=-3Q^{\musFlat{}}_{7} + Q^{\musFlat{}}_{21}
\end{equation*}

But since $Q^{\musFlat{}}_{7}$ and $Q^{\musFlat{}}_{21}$ are linearly independant they must individually be zero which implies that $Q$ itself is zero.

\begin{equation*}
Q^{\musFlat{}}_{7}=Q^{\musFlat{}}_{21}=Q^{\musFlat{}}=Q=0
\end{equation*}

Therefore, the only two-multivector field mapped to zero is the zero multivector field which means the map is injective. To prove surjectivity we need to demonstrate that any differential two-form $\beta$ may be expressed as the interior product of some two-multivector field $Q$ with the Cayley four-form. We start by decomposing the given two-form $\beta$ into irreducible Spin(7) representations.

\begin{equation*}
\beta=\beta_{7} + \beta_{21}
\end{equation*}

Using the representation equations we may express the components of $\beta$ in terms of the interior product.

\begin{align*}
& \beta_{7}=-\frac{1}{3}(\ast (\Psi \wedge \beta_{7})) = -\frac{1}{3} \beta^{\musSharp{}}_{7} \lrcorner \Psi \\
& \beta_{21}=\ast (\Psi \wedge \beta_{21}) = \beta^{\musSharp{}}_{21} \lrcorner \Psi 
\end{align*}

All we need to do from here to get an expression for $\beta$ is to sum the equations. 

\begin{equation*}
\beta=\big(-\frac{1}{3}\beta^{\musSharp{}}_{7}+\beta^{\musSharp{}}_{21}\big) \lrcorner \Psi
\end{equation*}

Therefore, any two-form may be expressed as the interior product of a two-multivector field with the Cayley form, and we know how to relate the corresponding multivector field to the original two-form. As any two-form may be expressed in this manner the map $\Psi_{2}$ is also surjective. \\

\begin{remark} All of these results concerning the maps $\Psi_{i}$ follow from Schur's Lemma. In particular, Schur's Lemma implies that the map $\Psi_{1}$ is an isomorphism onto its image $\Lambda^{3}_{7}(T^{*}M)$. Moreover, $\Psi_{3}$ is an isomorphism when restrited to $\Lambda^{3}_{8}(T^{*}M)$ and it maps any element of $\Lambda^{3}_{48}(T^{*}M)$ to zero.  We present the result in this way to make clear precisely how these isomorphisms works and to determine the constants which are not clear from Schur's Lemma on abstract grounds. 
\end{remark}

Now we discuss Hamiltonian multivector fields and differential forms corresponding to a torsion-free Spin(7)-Structure, which we refer to as \textit{Cayley multivector fields} and \textit{Cayley differential forms}. \\

We have seen that the map $\Psi_{3}$ is surjective which means that given an exact one-form df there exists at least one three-multivector field $Q$ which is Cayley and related to $df$ through the Cayley four-form. 

\begin{equation*}
Q \lrcorner \Psi = df
\end{equation*}

Similarly, we saw that the map $\Psi_{2}$ is an isomorphism so given an exact two-form $d\alpha$ there exists a unique two-multivector field Q which is Cayley corresponding to $d\alpha$.

\begin{equation*}
Q \lrcorner \Psi = d\alpha
\end{equation*}

However, all we know about the map $\Psi_{1}$ is that it is injective, so it is not immediately clear if nontrivial Cayley vector fields exist. It turns out that they do not.

\begin{theorem}
There do not exist any non-trivial Cayley vector fields on a closed Spin(7)-Manifold.
\end{theorem}

\textit{Proof.} As a first step we prove that the Lie derivative of the Cayley four-form with respect to a Cayley vector field must vanish. By definition for a Cayley vector field X there is an exact three-form $d\beta$ that is equal to the contraction of X with the Cayley form. 

\begin{equation*}
X \lrcorner \Psi = d\beta
\end{equation*}

From Cartan's Identity we may relate the Lie derivative to the interior product and exterior derivative.

\begin{equation*}
L_{X}\Psi=d(X \lrcorner \Psi) + X \lrcorner d\Psi
\end{equation*}

The first term vanishes because X is a Cayley vector field by assumption so $X\lrcorner\Psi$ is exact and the second term vanishes because the Spin(7)-Structure is assumed to be torsion-free so $\Psi$ is closed. \\

Since the Cayley four-form determines the metric, a Cayley vector field is a Killing vector field. As the Spin(7) metric is Ricci-Flat we may cite Bochner's theorem which then implies that the Cayley vector field is parallel \cite{JoyceBook}. If the holonomy is precisely Spin(7) we don't need to do anything else because the Cheeger-Gromoll Splitting Theorem implies that a compact Spin(7) holonomy manifold cannot have nontrivial parallel vector fields \cite{JoyceBook}. \\

However in the case that the holonomy is properly contained in Spin(7) we need to do a little more work. In this case the Cayley vector field is still parallel and this implies that the dual one-form $X^{\musFlat{}}$ is parallel. As the exterior derivative is the skew symmetrization of the covariant derivative since $X^{\musFlat{}}$ is parallel it is also closed. \\

In this case we utilize a projection identity relating vector fields and the Cayley four-form \cite{Karigiannis}.

\begin{equation*}
(X \lrcorner \Psi) \wedge \Psi = 7 \ast X^{\musFlat{}}
\end{equation*}

Wedging both sides of this equation by $X^{\musFlat{}}$ we get an expression for the pointwise norm of $X$.

\begin{equation*}
X^{\musFlat{}} \wedge (X \lrcorner \Psi) \wedge \Psi = 7 \| X \|^2 \textnormal{\textnormal{vol}}
\end{equation*}

We can then substitute in $d\beta$ for $X\lrcorner\Psi$ and rewrite the entire left-hand side as an exact eight-form since both $\Psi$ and $X^{\musFlat{}}$ are closed.

\begin{equation*}
d(X^{\musFlat{}} \wedge \beta \wedge \Psi) = 7 \| X \|^2 \textnormal{\textnormal{vol}}
\end{equation*}

By integrating both sides and using Stokes' Theorem we can see that $X=0$.  \qed \\

\begin{remark} A result similar to Theorem 4 was proved by Cho, Salur, and Todd for certain G$_{2}$-Structures starting from the identity which gives the G$_{2}$-metric in terms of the G$_{2}$ three-form. The Cayley four-form also determines a metric though it does so in a complicated enough manner that no direct analogue of the G$_{2}$ proof is available in the Spin(7) context. To prove our result we had to utilize results about the curvature of a torsion-free Spin(7) structure. Our method of proof also extends to the G$_{2}$ case, however, due to the curvature assumptions this method only applies to torsion-free G$_{2}$-Structures whereas the original proof only required that the G$_{2}$-Structure be closed. Even still there is no lower torsion class for a Spin(7)-Structure that would have the Cayley form be closed so we don't lose anything by using these assumptions on the curvature.
\end{remark}

We can apply a similar methodology to study Cayley two-multivector fields. We can be confident that they exist in this case as one will correspond to the exterior derivative of any one-form on the manifold. More generally, for any two-multivector field $Q$ on a compact manifold we may use the Hodge decomposition to write the contraction of the two-multivector field into the Cayley form uniquely as the sum of an exact, co-exact, and harmonic contribution.

\begin{equation*}
Q \lrcorner \Psi = d\alpha + d^{\ast}\beta + \gamma
\end{equation*}

From this perspective, the term $d^{\ast}\beta$ is what would prevent Q from being locally Cayley and $\gamma$ then prevents it from being Cayley. Note that if the second Betti number vanishes then there cannot be a harmonic contribution in the Hodge decomposition. There are examples of compact manifolds satisying this topological condtion \cite{Joyce1, Joyce2}. On such manifolds every locally Cayley two-multivector field is Cayley.\\

Now we assume we are given a Cayley two-multivector field $Q$ with corresponding Cayley form $d\alpha$. We may decompose $Q$ into irreducible Spin(7) representations.

\begin{equation*}
Q^{\musFlat{}}= Q^{\musFlat{}}_{7}+Q^{\musFlat{}}_{21}
\end{equation*}

We may use the precise form of the isomorphism $\Psi_{2}$ to express $d\alpha$ in terms of these components.

\begin{equation*}
-3  Q^{\musFlat{}}_{7} + Q^{\musFlat{}}_{21} = d\alpha
\end{equation*}

We may then differentiate this expression to find a constraint satisfied by all Cayley two-multivector fields.

\begin{equation*}
3 dQ^{\musFlat{}}_{7} = dQ^{\musFlat{}}_{21}
\end{equation*}

From the perspective of the Hodge decomposition this constraint is equivalent to $Q \lrcorner \Psi$ not having a co-exact contribution since, on a compact manifold, $d^{\ast}\beta$ vanishes if and only if $dd^{\ast}\beta$ vanishes. \\

We now consider some constraints on Cayley two-multivector fields. We first prove that a Cayley two-multivector field on a compact manifold cannot be the wedge product of two vector fields, in other words, it cannot be decomposable. 

\begin{proposition}
A Cayley two-multivector field on a closed Spin(7)-Manifold cannot be decomposable. 
\end{proposition}

\textit{Proof.} We start from the fundamental Spin(7) identity for vector fields $u$ and $v$ \cite{Karigiannis}.

\begin{equation*}
(u \lrcorner v \lrcorner \Psi) \wedge (u \lrcorner v \lrcorner \Psi) \wedge \Psi = -6 \| u^{\musFlat{}} \wedge v^{\musFlat{}} \|^2 \textnormal{\textnormal{vol}}
\end{equation*}

Denote the Cayley two-multivector field by Q and the corresponding Cayley one-form by $\alpha$ such that $Q\lrcorner\Psi$ = $d\alpha$. As $Q$ is decomposable we may write it as the wedge product of two vector fields $u$ and $v$.

\begin{equation*}
Q= u \wedge v
\end{equation*}

We can then plug in $d\alpha$ in the appropriate places in the expression for $\| u^{\musFlat{}} \wedge v^{\musFlat{}} \|^2$.

\begin{equation*}
d\alpha \wedge d\alpha \wedge \Psi = d(\alpha \wedge d\alpha \wedge \Psi) = -6  \| u^{\musFlat{}} \wedge v^{\musFlat{}} \|^2 \textnormal{\textnormal{vol}}
\end{equation*}

Applying Stokes' Theorem as before we can conclude that $u^{\musFlat{}} \wedge v^{\musFlat{}} = 0$. Since $\musFlat{}(u\wedge v) = u^{\musFlat{}} \wedge v^{\musFlat{}}$ is an isomorphism this implies that $u \wedge v=0$. \qed \\

This result  can be thought of as a restriction on differential forms on a Spin(7)-Manifold in the sense that a Spin(7)-Manifold does not admit exact two-forms whose preimage under the isomorphism $\Psi_{2}$ is decomposable. \\

We can now relate the norms of $Q^{\musFlat{}}_{7}$, $Q^{\musFlat{}}_{21}$, and $d\alpha$ on a compact manifold.

\begin{lemma}
For a two-multivector field $Q$ on a Spin(7)-Manifold Q$^{\musFlat{}}$ and the norms of its irreducible Spin(7) components are related through an expression involving the Cayley form. 

\begin{equation*}
(Q \lrcorner \Psi) \wedge (Q \lrcorner \Psi) \wedge\Psi = (-27 \| Q^{\musFlat{}}_{7} \|^{2} + \| Q^{\musFlat{}}_{21} \|^{2} ) \textnormal{\textnormal{vol}}
\end{equation*}
\end{lemma}

\textit{Proof.} We start by decomposing the two-form $Q\lrcorner\psi$ into irreducible Spin(7) representations. 

\begin{equation*}
Q \lrcorner \Psi = \ast (Q^{\musFlat{}} \wedge \Psi) = -3 Q^{\musFlat{}}_{7} + Q^{\musFlat{}}_{21} 
\end{equation*}

Then we may take the wedge product of this decomposition with $\psi$ and then with $Q\lrcorner\psi$. Using the equations defining the representations and orthogonality we arrive at the required result. 

\begin{align*}
& (Q \lrcorner \Psi) \wedge \Psi = 9 \ast Q^{\musFlat{}}_{7} + \ast Q^{\musFlat{}}_{21} \\
& (Q \lrcorner \Psi) \wedge (Q \lrcorner \Psi) \wedge \Psi = (-27 \| Q^{\musFlat{}}_{7} \|^{2} + \| Q^{\musFlat{}}_{21} \|^{2}) \textnormal{\textnormal{vol}} 
\end{align*}

\qed

\begin{theorem}
For a Cayley Two-Multivector Field $Q$ on a closed Spin(7)-Manifold the L$^{2}$ norms of the components of the irreducible Spin(7) representation components of $Q$ are proportional to each other up to a factor of 27. 
\begin{equation*}
 27 \| Q^{\musFlat{}}_{7} \|^{2}_{L^{2}} = \| Q^{\musFlat{}}_{21} \|^{2}_{L^{2}}
\end{equation*}
\end{theorem}

\textit{Proof.} As the two-multivector field $Q$ is Cayley there is a one-form $\alpha$ that is related to it in the usual way.

\begin{equation*}
Q \lrcorner \Psi = d\alpha
\end{equation*}

We can then plug this into the result from the above lemma and utilize Stokes' Theorem. 

\begin{align*}
& d\alpha \wedge d\alpha \wedge \Psi = (-27 \| Q^{\musFlat{}}_{7} \|^{2} + \| Q^{\musFlat{}}_{21} \|^{2}) \textnormal{\textnormal{vol}} \\
& d( \alpha \wedge d\alpha \wedge \Psi) = (-27 \| Q^{\musFlat{}}_{7} \|^{2} + \| Q^{\musFlat{}}_{21} \|^{2}) \textnormal{\textnormal{vol}} 
\end{align*}

\qed

An implication of this result is that a Cayley two-multivector field on a compact manifold is either identically zero or it has components in both of the irreducible Spin(7) representations on two-forms. Also in the case that the left-hand side vanishes without having to integrate, as would be the case for instance if $d\alpha$ was decomposable, the squared  norms of its irreducible Spin(7) representations would be related pointwise up to a factor of 27. \\

Since the isomorphism $\Psi_{2}$ maps $Q_{7}$ to $-3 Q^{\musFlat{}}_{7}$ and maps $Q_{21}$ to $Q^{\musFlat{}}_{21}$ we may infer from this result that the L$^{2}$ norms of the components of $d\alpha$ are related up to a factor of nine. 

\begin{equation*}
9\|(d\alpha)_{7}\|_{L^2}^2 = \|(d\alpha_{21})\|_{L^2}^2
\end{equation*}

We may also prove that if $Q$ is a Cayley two-multivector field then $Q^{\musFlat{}}$ cannot be closed. 

\begin{proposition} If Q is a Cayley two-multivector field on a closed manifold then $ Q^{\musFlat{}}$ cannot be closed. Moreover, neither $ Q^{\musFlat{}}_{7}$ nor $ Q^{\musFlat{}}_{21}$ can be closed.
\end{proposition}

\textit{Proof}. We may express the closed condition on $Q^{\musFlat{}}$ in terms of its Spin(7) decomposition.

\begin{equation*}
dQ^{\musFlat{}} = dQ^{\musFlat{}}_{7} + dQ^{\musFlat{}}_{21} = 0
\end{equation*}

However, we already saw that these derivatives satisfy a similar linear condition from $\Psi_{2}$.

\begin{equation*}
3 dQ^{\musFlat{}}_{7} = dQ^{\musFlat{}}_{21}
\end{equation*}

Combining these equations we see that $Q^{\musFlat{}}$ being closed implies both of its irreducible Spin(7) components are individually closed. Looking at the equations which define their representations and utilizing the fact that $\Psi$ is torsion-free it follows that not only are $Q^{\musFlat{}}_{7}$ and $Q^{\musFlat{}}_{21}$ closed but they are also co-closed. Neither of these forms can individually vanish by Theorem 7 and they are distinct since they are in different Spin(7) representations. Therefore we have two distinct harmonic forms related by an exact form. Such an expression contradicts the uniqueness of the harmonic representative of a cohomology class on a compact manifold. We may also see from these manipulations that $ Q^{\musFlat{}}_{7}$ being closed implies that $Q^{\musFlat{}}_{21}$ is closed and vice versa. So if either one of these is closed then $Q^{\musFlat{}}$ is closed and we just demonstrated that $ Q^{\musFlat{}}$ cannot be closed.\qed \\

Note also that as a corollary $Q^{\musFlat{}}$ cannot be the exterior derivative of a one-form. \\

Lastly we consider Cayley three-multivector fields. The interior product of such a multivector field $Q$ with the Cayley four-form gives the differential of a function on the manifold.

\begin{equation*}
Q \lrcorner \Psi =df
\end{equation*}

As the map $\Psi_{3}$ is surjective, for a given $df$ there is at least one $Q$ that satisfies this equation. More generally, for an arbitrary three-multivector field on a compact manifold we may decompose its contraction into the Cayley form according to the Hodge decomposition. If we assume full Spin(7) holonomy we get a simplification because in this case the first Betti number vanishes. As a result there cannot be a harmonic term in the Hodge decomposition.

\begin{equation*}
Q \lrcorner \Psi = df+d^{\ast}\beta
\end{equation*}

In this case every locally Cayley three-multivector field is Cayley and the term $d^{\ast}\beta$ is what prevents an arbitrary three-multivector field from being Cayley.\\

On a compact manifold $d^{\ast}\beta$ vanishing is equivalent to $dd^{\ast}\beta$ vanishing. We may get an expression for the latter by differentiating the contraction of Q into $\Psi$.

\begin{equation*}
d(Q \lrcorner \Psi) = dd^{\ast}\beta
\end{equation*}

To simplify this expresson it is useful to decompose the three-multivector field Q into irreducible Spin(7) representations.

\begin{equation*}
Q^{\musFlat{}} = Q^{\musFlat{}}_{8} + Q^{\musFlat{}}_{48}
\end{equation*}

 It follows from Schur's Lemma that a three-multivector field in $\Lambda^{3}_{48}(T^{*}M)$ contracted into $\Psi$ always vanishes so the resulting one-form is determined by the $\Lambda^{3}_{8}(T^{*}M)$ component of the three-multivector field. We may then rewrite the constraint on $Q$ in terms of its $\Lambda^{3}_{8}(T^{*}M)$ component.

\begin{equation*}
d \ast (Q^{\musFlat{}} \wedge \Psi) = d \ast (Q^{\musFlat{}}_{8} \wedge \Psi) = dd^{\ast}\beta = 0
\end{equation*}

This constraint can be shown to imply that $Q^{\musFlat{}}_{8}$ must be co-closed. We may actually prove an even stronger result which gives an expression for $Q^{\musFlat{}}_{8}$ in terms of the Cayley form and Cayley function $f$.

\begin{proposition} For a Cayley three-multivector field Q the component $Q^{\musFlat{}}_{8} $ is co-exact and may be expressed succinctly in terms of the function f and Cayley four-form $\Psi$. Moreover, the pointwise norms of the Cayley three-multivector field and Caley one-form are related up to a factor of seven. 
\end{proposition}

\textit{Proof.} We start with a general Spin(7) identity and substitute $df$ for $X^{\musFlat{}}$ \cite{Karigiannis}.

\begin{equation*}
df \wedge \Psi = 7 \ast Q^{\musFlat{}}_{8} 
\end{equation*}

Then we use the fact that $\Psi$ is torsion-free to simplify the wedge product of $df$ with $\Psi$.

\begin{equation*}
d(f\Psi) = 7 \ast Q^{\musFlat{}}_{8} 
\end{equation*}

Now we can utilize the fact the $\Psi$ is self-dual and we can then take the Hodge star of both sides to express $Q^{\musFlat{}}_{8}$ as a co-exact form.

\begin{equation*}
d^{\ast} (f\Psi) = 7Q^{\musFlat{}}_{8} 
\end{equation*} 

To get the norm result we need only take the norm of the original identity.

\begin{equation*}
df \wedge \Psi \wedge \ast(df \wedge \Psi) = 7 \| Q^{\musFlat{}}_{8}\|^2
\end{equation*}

We may simplify the left-hand side into an expression for the norm of $df$. \cite{Karigiannis}

\begin{equation*}
\| df \|^2 = 7 \| Q^{\musFlat{}}_{8}\|^2
\end{equation*}

\qed

We saw that for a Cayley two-multivector field on a compact manifold its metric dual cannot be closed. We may prove an analogous result for Cayley three-multivector fields. 

\begin{proposition} If Q is a Cayley three-multivector field on a closed manifold then neither $Q^{\musFlat{}}$ nor its component $Q^{\musFlat{}}_{8}$ can be closed. 
\end{proposition}

\textit{Proof}. We first write out an expression for the norm of the Cayley one-form. 

\begin{equation*}
Q^{\musFlat{}} \wedge (Q \lrcorner \Psi) \wedge \Psi = Q^{\musFlat{}}_{8} \wedge (Q \lrcorner \Psi) \wedge \Psi = 7 \| df \|^{2} \textnormal{vol}
\end{equation*}

The first equality holds because $Q^{\musFlat{}}_{48}$ vanishes when wedged against $\Psi$. If $Q$ is Cayley and either $Q^{\musFlat{}}$ or $Q^{\musFlat{}}_{8}$ are closed we may use Stokes theorem to show that $df$ vanishes. \qed

\begin{remark}
It can also be shown that $Q^{\musFlat{}}_{8}$ being closed is equivalent to the function f being harmonic. As a harmonic function on a compact manifold is constant, this would also imply that $df$ would necessarily vanish.
\end{remark}

Following from these results, there is generally a relationship between $df$ and $Q^{\musFlat{}}_{8}$. A priori, it may not be expected though that $Q^{\musFlat{}}$ itself would be constrained by $df$. After all, given a Cayley three-multivector field $Q$ with Cayley function $f$ we can add to $Q$ any element of $\Lambda^{3}_{48}$ and we will get a Cayley multivector field with the same Cayley function. However, the norm of $Q$ does turn out to be related to the norm of $df$ through the Spin(7) triple product operation. This triple product maps three vector fields $u,v$, and $w$ to a vector field $X$ utilizing the interior product against $\Psi$.

\begin{equation*}
X(u,v,w)^{\musFlat{}}=w \lrcorner v \lrcorner u \lrcorner \Psi
\end{equation*}

By extending the triple product operation by linearity we may define a vector field X(Q) given any three-multivector field. If the corresponding one-form is Cayley then $X(Q)=df$. We may now consider the relationship between $X(Q), Q^{\musFlat{}}_{8}, \textnormal{and}  Q^{\musFlat{}}_{48}$. This relationship is particularly simple when $Q$ is decomposable. 

\begin{proposition}
For a decomposable Cayley three-multivector field Q the norms of the irreducible Spin(7) components may be given in terms of constant multiples of the norm of $Q^{\musFlat{}}$.

\begin{align*}
& \|Q^{\musFlat{}}_{8} \|^{2} =1/7 \hspace{0.5mm} \|Q^{\musFlat} \|^{2} \\
& \|Q^{\musFlat{}}_{48} \|^{2} = 6/7 \hspace{0.5mm} \|Q^{\musFlat} \|^{2} 
\end{align*}
\end{proposition}

\textit{Proof.} We start with the relationship of the standard triple product norm and the norm of $Q=u\wedge v \wedge w$ \cite{Karigiannis}. Since we presume $Q$ is Cayley we let $df$ denote the one-form resulting from the triple product of $u,v, \textnormal{and} w$.

\begin{equation*}
\| df \|^{2} = \| Q^{\musFlat{}} \|^{2}
\end{equation*}

However, we already know how to relate $df$ to $Q^{\musFlat{}}_{8}$ and moreover, we know that the norm of $Q$ is equal to the Pythagorean sum of the norms of its irreducible Spin(7) components. The result then follows.\qed \\

From these relations, we can see that a decomposable three-multivector field has nonzero components in both irreducible Spin(7) representations on three-forms. If we also assume compactness then we know that there are points on the manifold where df must vanish. For a general Cayley multivector field $ Q^{\musFlat{}}_{8}$ must vanish at these points as well. However if we also assume that the Cayley three-multivector field is decomposable then these results on the norms imply that at the points where $df$ is zero $Q$ must vanish identically. So on a compact manifold a decomposable Cayley three-multivector field cannot be nonvanishing. \\

Cayley three-multivector fields in general seem to be pretty weak geometric objects due to the freedom to add any element of $Q^{\musFlat{}}_{48}$ to a Cayley three-multivector field while preserving the Cayley condition. Although adding the decomposability requirement constrains these structures considerably, given the flexibility of these structures in general it is likely that decomposable Cayley three-multivector fields exist, even on compact manifolds. What is certainly true is that if we remove the compactness assumption we may easily write down examples of these structures. To construct them we turn to our local model of $\mathbb{R}^{8}$ with the canonical Spin(7)-Structure $\Psi_{0}$. We then take $Q$ as the product of three coordinate vector fields.

\begin{equation*}
Q= \partial_{x_{0}} \wedge \partial _{x_{1}} \wedge \partial_{x_{2}}
\end{equation*}

Then we may compute $Q \lrcorner \Psi_{0}$. 

\begin{equation*}
Q \lrcorner \Psi_{0} = \ast (Q^{\musFlat{}} \wedge \Psi) = \ast (dx_{0124567}) =dx_{3}
\end{equation*}

As $dx_{4}$ is the derivative of a coordinate function $Q$ is a decomposable Cayley three-multivector field. Note that we cannot construct decomposable Cayley two-multivector fields so easily basically because not enough terms are killed off when wedging with the Cayley form in that case to get an exact form.  \\

%%%%%%%%%%%%%%%%%%%%%%%%%%%%%%%%%%%%%%%%%%%%%%%%%%%%%%%%%%%%%%%%%%

\subsection*{Acknowledgements} I would like to thank my advisor Dave Morrison for many helpful discussions and the Simons Collaboration on Special holonomy in Geometry, Analysis, and Physics for many interesting talks which inspired me to study the connections between exceptional holonomy and symplectic geometry.  I am also very appreciative to my friends and family, in particular Jessica Li, for being tolerant of the level of abstraction through which I approach the world. The research of the author is partially suported by Simons Foundation Award \#488629 (Morrison)

 \end{document}